\documentclass[12pt]{amsart}
\usepackage[margin=1in]{geometry}

\newtheorem{theorem}{Theorem}[section]
\newtheorem{lemma}[theorem]{Lemma}
\newtheorem{conjecture}[theorem]{Conjecture}

\newtheorem{corollary}[theorem]{Corollary}
\newtheorem{proposition}[theorem]{Proposition}

\theoremstyle{definition}

\newtheorem{example}[theorem]{Example}

\theoremstyle{remark}
\newtheorem{remark}[theorem]{Remark}

\begin{document}

\title[]
{Revitalized automatic proofs: \\
demonstrations}



\author{Tewodros Amdeberhan, David Callan, \\
Hideyuki Ohtsuka and
Roberto Tauraso}
\address{Department of Mathematics,
Tulane University, New Orleans, LA 70118, USA}
\email{tamdeber@tulane.edu}

\address{Department of Statistics, 
University of Wisconsin-Madison, Madison, WI 53706, USA}
\email{callan@stat.wisc.edu}

\address{Bunkyo University High School, 1191-7, 
Kami, Ageo-city, Saitama Pref., 362-0001, Japan}
\email{otsukahideyuki@gmail.com}

\address{Dipartimento di Matematica,
 Universita' di Roma "Tor Vergata", 
 Via della Ricerca Scientifica, 1, 
 00133 Roma, Italy}
\email{tauraso@mat.uniroma2.it}


\date{\today}


\begin{abstract}
We consider three problems from the recent issues of the American Mathematical Monthly involving different versions of Catalan triangle. Our main results offer generalizations of these identities and demonstrate automated proofs with additional twists, and on occasion we furnish a combinatorial proof. 
\end{abstract}

\maketitle

\newcommand{\ba}{\begin{eqnarray}}
\newcommand{\ea}{\end{eqnarray}}
\newcommand{\ift}{\int_{0}^{\infty}}
\newcommand{\nn}{\nonumber}
\newcommand{\no}{\noindent}
\newcommand{\lf}{\left\lfloor}
\newcommand{\rf}{\right\rfloor}
\newcommand{\realpart}{\mathop{\rm Re}\nolimits}
\newcommand{\imagpart}{\mathop{\rm Im}\nolimits}

\newcommand{\op}[1]{\ensuremath{\operatorname{#1}}}
\newcommand{\pFq}[5]{\ensuremath{{}_{#1}F_{#2} \left( \genfrac{}{}{0pt}{}{#3}
{#4} \bigg| {#5} \right)}}

\newtheorem{Definition}{\bf Definition}[section]
\newtheorem{Thm}[Definition]{\bf Theorem}
\newtheorem{Example}[Definition]{\bf Example}
\newtheorem{Lem}[Definition]{\bf Lemma}
\newtheorem{Cor}[Definition]{\bf Corollary}
\newtheorem{Prop}[Definition]{\bf Proposition}
\numberwithin{equation}{section}

\section{Introduction}

\smallskip
\noindent
Let's fix some nomenclature. The set of all integers is $\mathbb{Z}$, and the set of non-negative integers is $\mathbb{N}$. Denote the Catalan triangle by $B_{n,k}=\frac{k}n\binom{2n}{n-k}=\frac{k}n\binom{2n}{n+k}$, for $1\leq k\leq n$, and the all-familiar Catalan numbers $C_n=\frac1{n+1}\binom{2n}n$ correspond to $B_{n,1}$. On the other hand, $t_{2n-k,k}=\binom{2n}k-\binom{2n}{k-1}$ form yet another variation of the Catalan triangle and these numbers count lattice paths (N and E unit steps) from $(0,0)$ to $(2n-k,k)$ that may touch but stay below the line $y=x$. 

\smallskip
\noindent
\bf Convention. \rm Empty sums and empty products are evaluated to $0$ and $1$, respectively. Also that $\binom{n}k=0$ whenever $k<0$ or $k>n$.

\smallskip
\noindent
Let $Q_{a,b}:=\binom{a+b}a$. When considering a triple product of the numbers $B_{n,k}$, on occasion we find the following as a more handy reformulation
\begin{align} \label{GE-1}
\frac{abc\,Q_{a,b}Q_{b,c}Q_{c,a}}{Q_{a,a}Q_{b,b}Q_{c,c}}B_{a,k}B_{b,k}B_{c,k}=
 k^3\binom{a+b}{a+k}\binom{b+c}{b+k}\binom{c+a}{c+k}.
\end{align}

\smallskip
\noindent
The impetus for this paper comes from Problem $11844$ ~\cite{A}, Problem $11899$ ~\cite{B} and Problem $11916$ \cite{C} of the American Mathematical Monthly journal, plus the following identities that came up in our study:
\begin{align} \label{AMM1}
\binom{n+m}{2n}\sum_{k=0}^nk\binom{2n}{n+k}^2\binom{2m}{m+k}=\frac{n}2\binom{2m}{m+n}\binom{2n}n\sum_{j=0}^{m-1}\binom{n+j}n\binom{n+j}{n-1},
\end{align}
\begin{align} \label{AMM2}
\binom{n+m}m\sum_{k=0}^nk\binom{2n}{n+k}\binom{2m}{m+k}^2
=\frac{n}2\binom{2n}n\binom{2m}m\sum_{j=0}^{m-1}\binom{n+j}n\binom{m+j}{m-1}.
\end{align}
The purpose of our work here is to present certain generalizations and to provide {\em automatic proofs} as well as alternative techniques. Our demonstration of the Wilf-Zeilberger style of proof \cite{WZ} exhibit the power of this methodology, especially where we supplemented it with novel adjustments whenever a direct implementation lingers. 

\smallskip
\noindent
A class of $d$-fold binomial sums of the type
$$R(n):=\sum_{k_1,\dots,k_d}\prod_{i=1}^d\binom{2n}{n+k_i}\vert f(k_1,\dots,k_d)\vert$$
have been investigated by several authors, see for example ~\cite{BOOP} and references therein. One interpretation is this: $4^{-dn}R(n)$ is the expectation of $\vert f\vert$ if one starts at the origin and takes $2n$ random steps $\pm\frac12$ in each of the $d$ dimensions, thus arriving at the point $(k_1,\dots,k_d)\in\mathbb{Z}^d$ with probability
$$4^{-dn}\prod_{i=1}^d\binom{2n}{n+k_i}.$$

\smallskip
\noindent
The organization of the paper is as follows. In Section 2, Problems 11844, 11916 and some generalized identities are proved. Section 3 resolves Problem 11899 and highlights a combinatorial proof together with $q$-analogue of related identities. Finally, in Section 4, we conclude with further generalizations and some open problems for the reader.

\section{The first set of main results}

\bigskip
\noindent
Our first result proves Problem $11844$ of the Monthly ~\cite{A} as mentioned in the Introduction.

\begin{lemma} \label{GT3} For non-negative integers $m\geq n$, we have
\begin{align} \label{GE3}
\sum_{k=0}^n(m-2k)\binom{m}k^3=(m-n)\binom{m}n\sum_{j=0}^{m-n-1}\binom{n+j}n\binom{n+j}{m-n-1}. \end{align}
\end{lemma}
\begin{proof} We apply the method of Wilf-Zeilberger ~\cite{WZ}. This techniques works, in the present case, after multiplying \eqref{GE3} through with $(-1)^m$. Denote the resulting summand on the LHS of \eqref{GE3} by $F_1(m,k)$ and its sum by $f_1(m):=\sum_{k=0}^nF_1(m,k)$. Now, introduce the companion function 
$$G_1(m,k):=-F_1(m,k)\cdot \frac{(2m-k+2)k^3}{(m-2k)(m-k+1)^3}$$
and check that $F_1(m+1,k)-F_1(m,k)=G_1(m,k+1)-G_1(m,k)$. Telescoping gives
\begin{align*} 
f_1(m+1)-f_1(m)&=\sum_{k=0}^nF_1(m+1,k)-\sum_{k=0}^nF_1(m,k)
=\sum_{k=0}^n[G_1(m,k+1)-G_1(m,k)] \\
&=G_1(m,n+1)-0=(-1)^{m+1}\binom{m}n^3(2m-n+1).
\end{align*}
Let $F_2(m,j)$ be the summand on the RHS of \eqref{GE3} and its sum $f_2(m):=\sum_{j=0}^{m-n-1}F_2(m,j)$. Introduce 
$$G_2(m,j):=F_2(m,j)\cdot\frac{j(m-2n-j-1)}{(m-n)^2}$$
and check that $F_2(m+1,j)-F_2(m,j)=G_2(m,j+1)-G_2(m,j)$. Summing $0\leq j\leq m-n$ and telescoping, we arrive at
\begin{align*} f_2(m+1)-f_2(m)&=\sum_{j=0}^{m-n}F_2(m+1,j)-\sum_{j=0}^{m-n}F_2(m,j)+F_2(m,m-n) \\
&=\sum_{j=0}^{m-n}[G_2(m,j+1)-G_2(m,j)]+F_2(m,m-n) \\
&=G_2(m,m-n+1)-0+F_2(m,m-n)\\
&=(-1)^{m+1}\binom{m}n^3(2m-n+1). 
\end{align*}
The final step is settled with $f_1(0)=f_2(0)=0$ (if $m=0$, so is $n=0$).
\end{proof}

\smallskip
\noindent
\begin{theorem} \label{GT0} For nonnegative integers $r, s$ and $m\geq n$, we have
\begin{align} \label{GE4} 
\sum_{k=0}^n\frac{(m-2k)\binom{m+r+s}{m,r,s}\binom{m}k\binom{m+2r}{k+r}\binom{m+2s}{k+s}}
{\binom{m+2r}{m+r}\binom{m+2s}{m+s}\binom{m+s}{n+s}}
=(m-n)\sum_{j=0}^{m-n+r-1}\binom{n+j}n\binom{n+j+s}{m-n+s-1}.
\end{align}
\end{theorem}
\begin{proof} Again we use the W-Z method. Multiply through equation \eqref{GE4} by $\binom{m+s}{n+s}$ and denote the summand on the new LHS of \eqref{GE4} by $F_1(r,k)$ and its sum by $f_1(r):=\sum_{k=0}^nF_1(r,k)$. Now, introduce the companion function 
$$G_1(r,k):=F_1(r,k)\cdot \frac{k(s+k)}{(m-2k)(m+r-k+1)}$$
and (routinely) check that $F_1(r+1,k)-F_1(r,k)=G_1(r,k+1)-G_1(r,k)$. Telescoping gives
\begin{align*} 
f_1(r+1)-f_1(r)&=\sum_{k=0}^nF_1(r+1,k)-\sum_{k=0}^nF_1(r,k)
=\sum_{k=0}^n[G_1(r,k+1)-G_1(r,k)] \\
&=G_1(r,n+1)-0=(m-n)\binom{m+s}{n+s}\binom{m+r}n\binom{m+r+s}{m-n+s-1}. 
\end{align*}
Denoting the entire sum on the RHS of \eqref{GE4} by $f_2(r)$, it is straightforward to see that
$$f_2(r+1)-f_2(r)=(m-n)\binom{m+s}{n+s}\binom{m+r}n\binom{m+r+s}{m-n+s-1}.$$
It remains to verify the initial condition $f_1(0)=f_2(0)$; that is,
\begin{align} \label{GE5}
\sum_{k=0}^n\frac{(m-2k)\binom{m+s}m\binom{m}k^2\binom{m+2s}{k+s}}{\binom{m+2s}{m+s}}
=(m-n)\binom{m+s}{n+s}\sum_{j=0}^{m-n-1}\binom{n+j}n\binom{n+j+s}{m-n+s-1}.
\end{align}
Denote the summand on the LHS of \eqref{GE5} by $F_2(s,k)$ and its sum by $f_2(s):=\sum_{k=0}^nF_2(s,k)$. Now, introduce the companion function 
$$G_2(s,k):=F_2(s,k)\cdot \frac{k^2}{(m-2k)(m+s-k+1)}$$
and (routinely) check that $F_2(s+1,k)-F_2(s,k)=G_2(s,k+1)-G_2(s,k)$. Telescoping gives
\begin{align*} 
f_2(s+1)-f_2(s)&=\sum_{k=0}^nF_2(s+1,k)-\sum_{k=0}^nF_2(s,k)
=\sum_{k=0}^n[G_2(s,k+1)-G_2(s,k)] \\
&=G_2(s,n+1)-0=(m-n)\binom{m+s}{n+s+1}\binom{m+s}n\binom{m}n. 
\end{align*}
Let $F_3(s,j)$ be the summand on the RHS of \eqref{GE5} and its sum $f_3(s):=\sum_{j=0}^{m-n-1}F_3(s,j)$. Introduce 
$$G_3(s,j):=F_3(s,j)\cdot\frac{j(2n-m+j+1)}{(n+s+1)(m-n+s)}$$
and check that $F_3(s+1,j)-F_3(s,j)=G_3(s,j+1)-G_3(s,j)$. Summing and telescoping, we get
\begin{align*} f_3(s+1)-f_3(s)&=\sum_{j=0}^{m-n-1}F_3(s+1,j)-\sum_{j=0}^{m-n-1}F_3(s,j)
=\sum_{j=0}^{m-n-1}[G_3(s,j+1)-G_3(s,j)] \\
&=G_3(s,m-n)-0=(m-n)\binom{m+s}{n+s+1}\binom{m+s}n\binom{m}n. 
\end{align*}
The initial condition $f_2(0)=f_3(0)$ is precisely the content of Lemma \ref{GT3}.
\end{proof}

\smallskip
\noindent
The next statement covers Problem $11916$ ~\cite{C} as an immediate application of Theorem \ref{GT0}.
 
\smallskip
\noindent
\begin{corollary}  \label{GT00} Let $a, b$ and $c$ be non-negative integers. Then, the function
$$U(a,b,c):=a\binom{a+b}a\sum_{j=0}^{c-1}\binom{a+j}a\binom{b+j}{b-1}$$
is symmetric, i.e. $U(\sigma(a),\sigma(b),\sigma(c))=U(a,b,c)$ for any $\sigma$ in the symmetric group $\mathfrak{S}_3$. 
\end{corollary}
\begin{proof} If $n=a, m=2a, r=b-a, s=c-a$, the left-hand side of Theorem \ref{GT0} turns into
\begin{align*} 
LHS&=\frac{\binom{b+c}{2a,b-a,c-a}}{\binom{2b}{b+a}\binom{2c}{c+a}}
\sum_{k=0}^a(2a-2k)\binom{2a}k\binom{2b}{k+b-a}\binom{2c}{k+c-a} \\
&=2\frac{(a+b)!(b+c)!(c+a)!}{(2a)!(2b)!(2c)!}\sum_{k=0}^ak\binom{2a}{a-k}\binom{2b}{b-k}\binom{2c}{c-k} \\
&=\frac{2Q_{a,b}Q_{b,c}Q_{c,a}}{Q_{a,a}Q_{b,b}Q_{c,c}}
\sum_{k=0}^ak\binom{2a}{a+k}\binom{2b}{b+k}\binom{2c}{c+k}
\end{align*}
and the right-hand side simplifies to
$$RHS =a\binom{a+c}c\sum_{j=0}^{b-1}\binom{a+j}a\binom{c+j}{c-1}=aQ_{c,a}\sum_{j=0}^{b-1}\binom{a+j}a\binom{c+j}{c-1}.$$
Therefore, we obtain
\begin{align} \label{GE0}
\frac{Q_{a,b}Q_{b,c}Q_{c,a}}{Q_{a,a}Q_{b,b}Q_{c,c}}
\sum_{k=0}^ak\binom{2a}{a+k}\binom{2b}{b+k}\binom{2c}{c+k}
=\frac{aQ_{c,a}}2\sum_{j=0}^{b-1}\binom{a+j}a\binom{c+j}{c-1}. \end{align}
The following apparently symmetry
$$\sum_{k=0}^ak\binom{2a}{a+k}\binom{2b}{b+k}\binom{2c}{c+k} 
=\sum_{k=0}^{\min\{a,b,c\}}k\binom{2a}{a+k}\binom{2b}{b+k}\binom{2c}{c+k}.$$ 
implies that the LHS of the identity in \eqref{GE0} has to be symmetric. The assertion follows from the symmetry inherited by the RHS of the same equation \eqref{GE0}.
\end{proof}

\smallskip
\noindent
\begin{example} In equation \eqref{GE0}, the special case $a=n, b=c=m$ becomes \eqref{AMM1} while $a=b=n, c=m$ recovers \eqref{AMM2}.
\end{example}

\smallskip
\noindent
\begin{corollary} \label{GT1} Preserve notations from Cor. \ref{GT00}. For $a, b, c\in\mathbb{N}$ and any $\sigma\in \mathfrak{S}_3$, we have
\begin{align} \label{GE1}
\sum_{k=0}^ak\binom{a+b}{a+k}\binom{b+c}{b+k}\binom{c+a}{c+k} 
=\frac{\sigma(a)Q_{\sigma(a),\sigma(b)}}2\sum_{j=0}^{\sigma(c)-1}
\binom{\sigma(a)+j}{\sigma(a)}\binom{\sigma(b)+j}{\sigma(b)-1}. \end{align}
\end{corollary}
\begin{proof} First, employ an algebraic manipulation on \eqref{GE0} similar to equation \eqref{GE-1}. Now apply the identity in \eqref{GE0} and the statement of Corollary \ref{GT00}.
\end{proof}

\noindent
For non-negative integers $x, y, z$, write the elementary symmetric functions 
$$e_1(x,y,z)=x+y+z, \qquad e_2(x,y,z)=xy+yz+zx \qquad \mbox{ and } \qquad e_3(x,y,z)=xyz.$$

\noindent
\begin{theorem} \label{GT2} For non-negative integers $a, b$ and $c$, we have
\begin{align} \label{GE2} 
\sum_{k=0}^ak^3\binom{a+b}{a+k}\binom{b+c}{b+k}\binom{c+a}{c+k}
=\frac{b^2c^2 Q_{b,c}}2\sum_{j=0}^{a-1}\frac{e_2(a,b,c)\binom{b+j}b\binom{c+j}c}{e_2(j,b,c)\cdot e_2(j+1,b,c)}.
\end{align}
\end{theorem}
\begin{proof} Once again use the W-Z method. First, divide through by $e_2(a,b,c)$ to denote the summand on the LHS of \eqref{GE2} by $F_1(a,k)$ and its sum by $f_1(a):=\sum_{k=0}^aF_1(a,k)$. 
Now, introduce the companion function 
$$G_1(a,k):=-F_1(a,k)\cdot 
\frac{((e_2+b+c)k^2-(e_2+b+c)k+abc+bc)(b+k)(c+k)}
{2k^3(a+1-k)\cdot(e_2+b+c)}$$
and (routinely) check that $F_1(a+1,k)-F_1(a,k)=G_1(a,k+1)-G_1(a,k)$; where we write $e_2$ for $e_2(a,b,c)$. Keeping in mind that $F_1(a,a+1)=0$ and telescoping gives
\begin{align*} 
f_1(a+1)-f_1(a)&=\sum_{k=0}^{a+1}F_1(a+1,k)-\sum_{k=0}^{a+1}F_1(a,k)
=\sum_{k=0}^{a+1}[G_1(a,k+1)-G_1(a,k)] \\ 
&=G_1(a,a+2)-G_1(a,0)=0-G_1(a,0)
=\frac{b^2c^2Q_{a,b}Q_{b,c}Q_{c,a}}{2e_2(e_2+b+c)}.
\end{align*}
This difference formula for $f_1(a+1)-f_1(a)$ leads to
$$f_1(a)=\frac{b^2c^2Q_{b,c}}2\cdot\sum_{j=0}^{a-1}\frac{\binom{b+j}a\binom{c+j}c}
{(jb+bc+cj)\cdot(jb+bc+cj+b+c)}$$
which is the required conclusion.
\end{proof}

\noindent
\begin{remark} In ~\cite{MOR}, Miana, Ohtsuka and Romero obtained two identities for the sum $\sum_{k=0}^nB_{n,k}^3$. 
From Theorem \ref{GT2} and \eqref{GE-1}, we 
obtain the identity for the sum $\sum_{k=0}^a B_{a,k}B_{b,k}B_{c,k}$.
\end{remark}

\smallskip
\noindent
\begin{remark} Corollary \ref{GT1} and Theorem \ref{GT2} exhibit formulas for $\sum_k k(\cdots)$ and $\sum_k k^3(\cdots)$. It appears that similar (albeit complicated) results are possible for sums of the type $\sum_k k^p(\dots)$ whenever $p$ is an odd positive integer (but not when $p$ is even).
\end{remark}

\smallskip
\noindent
We can offer a $4$-parameter generalization of Theorem \ref{GT1} and Theorem \ref{GT2}.

\smallskip
\noindent
\begin{theorem} \label{GT3.4} For non-negative integers $a, b, c $ and $d$, we have
$$\sum_{k=0}^a k\binom{a+b}{a+k}\binom{b+c}{b+k}\binom{c+d}{c+k}\binom{d+a}{d+k}
=\frac{bQ_{b,c}Q_{c,d}Q_{b+c+d,a}}{2Q_{a,c}}\sum_{j=0}^{a-1}\frac{Q_{b,j}Q_{c-1,j+1}Q_{d-1,j+1}}{Q_{b+c+d,j+1}}.$$
\end{theorem}
\begin{proof} Analogous to the preceding arguments. 
\end{proof}

\begin{remark} It is interesting to compare our results against Corollary 4.1 of ~\cite{GZ}. Although these are similar, there are differences: in our case the RHSs are less involved while those of ~\cite{GZ} are more general. See also Corollary 4.2 and Theorem 4.3 of ~\cite{MOR}. The examples below are devoted to explore some specifics.
\end{remark}

\begin{example} Set $a=b=c=n$ in Theorem \ref{GT1}. The outcome is
$$\sum_{k=0}^nk\binom{2n}{n+k}^3=\frac12\binom{2n}n\sum_{j=0}^nj\binom{n+j-1}{n-1}^2.$$
\end{example}

\begin{example} Set $a=b=c=n$ in Theorem \ref{GT2}. The outcome is
$$\sum_{k=0}^nk^3\binom{2n}{n+k}^3
=\frac12\sum_{j=0}^{n-1}\frac{3n^4\binom{2n}n\binom{n+j}n^2}{(n+2j)(n+2j+2)}.$$
\end{example}

\begin{example} Set $a=b=c=d=n$ in Theorem \ref{GT3.4}. The outcome is
$$\sum_{k=0}^nk\binom{2n}{n+k}^4
=\frac12\binom{4n}n\binom{2n}n
\sum_{j=0}^nj\binom{n+j-1}{n-1}^3\binom{3n+j}{3n}^{-1}.$$
\end{example}

\section{The second set of main results}

\smallskip
\noindent
We start with a $q$-identity and its ordinary counterpart will allow us to prove one of the Monthly problems which was alluded to in the Introduction. Along the way, we encounter the Catalan triangle $t_{2n-k,k}=\binom{2n}k-\binom{2n}{k-1}$ which we also write as  $t_{n+k,k}=\binom{2n}{n-k}-\binom{2n}{n-k-1}$.  Let's recall some notations. The $q$-analogue of the integer $n$ is given by $[n]_q:=\frac{1-q^n}{1-q}$, the factorial by $[n]_q!=\prod_{i=1}^n\frac{1-q^i}{1-q}$ and the binomial coefficients by
$$\binom{n}k_q=\frac{[n]_q!}{[k]_q![n-k]_q!}.$$

\begin{lemma} \label{GT4} For a free parameter $q$ and a positive integer $n$, we have 
\begin{align*} 
\sum_{k=0}^n\binom{2n+1}{n-k}_q\left[\binom{2n}{n-k}_q-\binom{2n}{n-k-1}_q\right]q^{k(k+1)}
=q^n\binom{2n}n_q^2.
\end{align*}
\end{lemma}
\begin{proof} Let $G(n,k)=\binom{2n}{n+k}_q^2\cdot q^{n+k^2}$. Now, check that 
\begin{align*}
\binom{2n+1}{n-k}_q\left[\binom{2n}{n-k}_q-\binom{2n}{n-k-1}_q\right]q^{k(k+1)}
&=\binom{2n+1}{n-k}_q^2\frac{1-q^{2k+1}}{1-q^{2n+1}}q^{n+k^2} \\
&=G(n,k)-G(n,k+1)
\end{align*}
and then sum over $k=0$ through $k=n$ to obtain $G(n,0)=q^n\binom{2n}n_q^2$.
\end{proof}

\smallskip
\noindent
We now demonstrate  a combinatorial argument for the special case $q=1$ of Lemma \ref{GT4}.
\begin{lemma} \label{GT5} For non-negative integers $n$, we have
\begin{align} \label{GE6}
\sum_{k=0}^n\binom{2n+1}{n-k}\left[\binom{2n}{n-k}-\binom{2n}{n-k-1}\right]=\binom{2n}n^2. \end{align}
\end{lemma} 
\begin{proof} The first factor in the summand on the left side of \eqref{GE6} counts paths of $2n+1$ steps, consisting of upsteps $(1,1)$ or downsteps $(1,-1)$, that start at the origin and end at height $2k+1$. The second factor is the generalized Catalan number that counts \it nonnegative \rm (i.e., first quadrant) paths of $2n$ up/down steps that end at height $2k$. By concatenating the first path and the reverse of the second, we see that the left side counts the set $X_n$ of paths of $2n+1$ upsteps and $2n$ downsteps that avoid the $x$-axis for $x>2n$, i.e. avoid $(2n+2,0),(2n+4,0),\dots,(4n,0)$. 

\smallskip
\noindent
Now $\binom{2n}{n}$ is the number of \it balanced \rm paths of length $2n$ (i.e., $n$ upsteps and $n$ downsteps), but it is also the number of nonnegative $2n$-paths and, for $n\geq1$, twice the number 
of positive (= nonnegative, no-return) $2n$-paths (see \cite{DC}, for example). So, the right side of \eqref{GE6} counts the set $Y_n$ of pairs $(P,Q)$ of nonnegative $2n$-paths. 
Here is a bijection $\phi$ from $X_n$ to $Y_n$. 
A path $P \in X_n$ ends at height 1 and so its last upstep from the $x$-axis splits it into $P=BUD$ where $B$ is a balanced path and $D$ is a dyck path of length $\ge 2n$ since $P$ avoids the $x$-axis for $x>2n$. 
Write $D$ as $QR$ where $R$ is of length $2n$. 

If $B$ is empty, set $\phi(P)=(Q$, Reverse($R$)), a pair of nonnegative $2n$-paths ending at the same height. If $B$ is nonempty, then by the above remarks it is equivalent to a bicolored positive path $S$ of the same length, say colored red or blue. If red, set $\phi(P)=(Q\,S$,  Reverse($R$))$\,\in Y_n$ with the first path ending strictly higher than the second.  If blue, set $\phi(P)=($Reverse($R$), $Q\,S$)$\,\in Y_n$ with the first path ending strictly lower than the second. It is easy to check that $\phi$ is a bijection from $X_n$ to $Y_n$.
\end{proof}

\smallskip
\noindent
As an application, we present a proof for Problem $11899$ as advertised in the Introduction.

\smallskip
\noindent
\begin{corollary} For non-negative positive integer $n$, we have
$$\sum_{k=0}^n\binom{2n}k\binom{2n+1}k+\sum_{k=n+1}^{2n+1}\binom{2n}{k-1}\binom{2n+1}k=\binom{4n+1}{2n}+\binom{2n}n^2.$$
\end{corollary}
\begin{proof} Start by writing
\begin{align*}
&A_1:=\sum_{k=0}^n\binom{2n}k\binom{2n+1}k, \qquad \qquad
A_2:=\sum_{k=n+1}^{2n+1}\binom{2n}k\binom{2n+1}k, \\ 
&\tilde{A}_1:=\sum_{k=n+1}^{2n+1}\binom{2n}{k-1}\binom{2n+1}k, 
\qquad \tilde{A}_2:=\sum_{k=0}^n\binom{2n}{k-1}\binom{2n+1}k. \end{align*}
Re-indexing gives $A_1=\tilde{A}_1$ and $A_2=\tilde{A}_2$. The required identity is $A_1+\tilde{A}_1=2A_1=\binom{4n+1}{2n}+\binom{2n}n^2$. In view of the Vandermonde-Chu identity $A_1+A_2=\binom{4n+1}{2n}$, it suffices to prove that
$A_1-A_2=A_1-\tilde{A}_2=\binom{2n}n^2$. That is, 
$$\sum_{k=0}^n\binom{2n+1}kt_{2n-k,k}
=\sum_{k=0}^n\binom{2n+1}{n-k}\left[\binom{2n}{n-k}-\binom{2n}{n-k-1}\right]=\binom{2n}n^2$$
which is exactly what Lemma \ref{GT5} is about. However, here is yet another verification:
if we let $G(n,k)=\binom{2n}{n+k}^2$ then it is routine to check that 
$$\binom{2n+1}{n-k}\left[\binom{2n}{n-k}-\binom{2n}{n-k-1}\right]
=\binom{2n+1}{n-k}^2\frac{2k+1}{2n+1}=G(n,k)-G(n,k+1).$$ 
Obviously then 
$$\sum_{k=0}^n[G(n,k)-G(n,k+1)]=G(n,0)-G(n,n+1)=G(n,0)=\binom{2n}n^2.$$
The proof follows.
\end{proof}

\section{Concluding Remarks}

\smallskip
\noindent
Finally, we list binomial identities with extra parameters similar to those from the preceding sections, however their proofs are left to the interested reader because we wish to limit unduly replication of our techniques. We also include some open problems.

\smallskip
\noindent
The first result generalizes Corollary \ref{GT1}.

\begin{proposition} For non-negative integers $a, b, c$ and an integer $r$, we have
\begin{align*}
\sum_{k=1}^{a+r}(2k-r)\binom{a+b+r}{a+k}\binom{b+c+r}{b+k}\binom{c+a+r}{c+k} 
=(a+r)Q_{a+r,b}\sum_{j=0}^{c+r-1}\binom{a+j}a\binom{b+j}{b+r-1}.
\end{align*}
\end{proposition}

\smallskip
\noindent
Next, we state certain natural $q$-analogues of Corollary \ref{GT1} and Corollary \ref{GT00}.

\smallskip
\noindent
\begin{theorem} For non-negative integers $a, b$ and $c$, we have
\begin{align*}
\sum_{k=0}^a\frac{(1-q^{2k})q^{2k^2-k-1}}{1-q^a}\binom{a+b}{a+k}_q\binom{b+c}{b+k}_q\binom{c+a}{c+k}_q
=\binom{a+b}a_q\sum_{j=0}^{c-1}q^j\binom{a+j}a_q\binom{b+j}{b-1}_q. 
\end{align*}
\end{theorem}

\begin{corollary} Let $a, b$ and $c$ be non-negative integers. Then, the function
\begin{align*} U_q(a,b,c):=\frac{1-q^a}{1-q}\binom{a+b}a_q\sum_{j=0}^{c-1}\binom{a+j}a_q\binom{b+j}{b-1}_q
\end{align*}
is symmetric, i.e. $U_q(\sigma(a),\sigma(b),\sigma(c))=U_q(a,b,c)$ for any $\sigma$ in the symmetric groups $\mathfrak{S}_3$.
\end{corollary}

\noindent 
Let's consider the family of sums
$$S_r(a,b,c):=\sum_{k=0}^a k^{2r+1}\binom{a+b}{a+k}\binom{b+c}{b+k}\binom{c+a}{c+k}.$$
It follows that
\begin{align*}
(a^2-k^{2})\binom{a+b}{a+k}\binom{c+a}{c+k}
&=(a+k)(a-k)\binom{a+b}{a+k}\binom{c+a}{a-k}\\
&=(a+b)(a+c)\binom{a-1+b}{a-1+k}\binom{c+a-1}{a-1-k}\\
&=(a+b)(a+c)\binom{a-1+b}{a-1+k}\binom{c+a-1}{c+k}
\end{align*}
which in turn implies, after replacing $k^{2r+1}=k^{2r-1}k^2=k^{2r-1}[a^2-(a^2-k^2)]$, that
$$S_r(a,b,c)=a^2\cdot S_{r-1}(a,b,c)-(a+b)(a+c)\cdot S_{r-1}(a-1,b,c).$$

\noindent
\bf Problem. \it  Introduce the operators on symmetric functions $f=f(a,b,c)$ of $3$-variables by 
$$\mathcal{L}\cdot f=[(a+b)(a+c)E - a^2 I]f$$
where $E\cdot f(a,b,c)=f(a-1,b,c)$ and $I\cdot f(a,b,c)=f(a,b,c)$ is the identity map. As a question of independent interest show that the iterates $\mathcal{L}^n\cdot 1$ always yield in symmetric polynomials in $\mathbb{Z}[a,b,c]$, for any integer $n\geq1$. \rm

\smallskip
\noindent
\bf Postscript. \rm Matthew Hongye Xie of Nankai University informed the authors, in private communication, that he has found a proof for this problem.

\smallskip
\noindent
\begin{conjecture} \label{conj} For each $r\in\mathbb{Z}^+$, there exist symmetric polynomials $f_r, g_r \in \mathbb{Z}[a,b,c]$ such that
\begin{align*}
S_r(a,b,c)=\frac{b^2c^2 f_r(a,b,c)\binom{b+c}{b}}{2}\sum_{j=0}^{a-1}\frac{\binom{b+j}b\binom{c+j}c\cdot g_r(j+1,b,c)}{f_r(j,b,c)\cdot f_r(j+1,b,c)}.
\end{align*}
The functions $f_r$ satisfy the recurrence, 
$$f_r(a,b,c)=\mathcal{L}\cdot f_{r-1}(a,b,c)$$
with $f_0(a,b,c)=1$.
\end{conjecture}

\begin{table}[h!]
  \centering 
  \label{tab:table1}  
  \begin{tabular}{|c|c|c|}
    \hline
    $r$ & $f_r$ & $g_r$\\
    \hline
    0 & 1 & $1/e_3$\\
    1 & $e_2$ & $1$\\
    2 & $e_2^2-e_1 e_2+e_3$ & $2e_3 -e_2$\\
    3 & $e_2^3+3e_3e_2-3e_2^2e_1-2e_3e_1+2e_2e_1^2+e_3-e_2e_1$      & $6e_3^2-8e_3e_2+3e_2^2+e_3-e_2e_1$\\
     \hline
  \end{tabular}
    \vspace{2mm}
   \caption{The first few polynomials in support of Conjecture \ref{conj}}
\end{table}

\end{document}